\documentclass[a4paper,11pt]{article}
\usepackage{amsmath,newtxtext,newtxmath,bm}
\usepackage[top=25mm,bottom=35mm,left=20mm,right=20mm,truedimen]{geometry}

\usepackage{amsthm}
\theoremstyle{definition}
\newtheorem{lem}{Lemma}[section]
\newtheorem{prop}{Proposition}[section]
\newtheorem*{pf}{Proof}
\usepackage{mathtools}
\mathtoolsset{showonlyrefs=true}
\usepackage{algorithm,algpseudocode}
\usepackage{hyperref}
\usepackage{graphicx}
\newcommand{\arsinh}{\mathrm{asinh}}
\newcommand{\artanh}{\mathrm{atanh}}
\newcommand{\sech}{\mathrm{sech}}

\title{Algorithms for the computation of the matrix logarithm based on the double exponential formula}
\author{
  Fuminori Tatsuoka\thanks{Department of Applied Physics, Graduate School of Engineering, Nagoya University, Furo-cho, Chikusa-ku, Nagoya 464-8603, Japan, \url{f-tatsuoka@na.nuap.nagoya-u.ac.jp}},~
  Tomohiro Sogabe\thanks{Department of Applied Physics, Graduate School of Engineering, Nagoya University, Furo-cho, Chikusa-ku, Nagoya 464-8603, Japan,~ \url{sogabe@na.nuap.nagoya-u.ac.jp}},~
  Yuto Miyatake\thanks{Cybermedia Center, Osaka University, 1-32 Machikaneyama, Toyonaka, Osaka 560-0043, Japan,~ \url{miyatake@cas.cmc.osaka-u.ac.jp}},~
  Shao-Liang Zhang\thanks{Department of Applied Physics, Graduate School of Engineering, Nagoya University, Furo-cho, Chikusa-ku, Nagoya 464-8603, Japan,~ \url{zhang@na.nuap.nagoya-u.ac.jp}}
}

\begin{document}
\maketitle
\begin{abstract}
  We consider the computation of the matrix logarithm by using numerical quadrature.
  The efficiency of numerical quadrature depends on the integrand and the choice of quadrature formula.
  The Gauss--Legendre quadrature has been conventionally employed;
  however, the convergence could be slow for ill-conditioned matrices.
  This effect may stem from the rapid change of the integrand values.
  To avoid such situations, we focus on the double exponential formula, which has been developed to address integrands with endpoint singularity.
  In order to utilize the double exponential formula, we must determine a suitable finite integration interval, which provides the required accuracy and efficiency.
  In this paper, we present a method for selecting a suitable finite interval based on an error analysis as well as two algorithms, and one of these algorithms is an adaptive quadrature algorithm.
\end{abstract}

\section{Introduction}
A logarithm of $A\in\mathbb{R}^{n\times n}$ is defined as any matrix $X$ such that
\begin{align}\label{eq:def_mat_log}
  \exp(X)=A,\end{align} where
$\exp(X):=I+X+\frac{1}{2!}X^2+\frac{1}{3!}X^3+\cdots$
\cite[p.~269]{higham2008functions}.
If all eigenvalues of $A$ lie in the set $\mathbb{C}\setminus (-\infty,0]$, there is a unique logarithm of $A$ whose eigenvalues all lie in the strip $\{z\in\mathbb{C}:-\pi<\mathrm{Im}(z)< \pi\}$ \cite[Thm. 1.31]{higham2008functions}.
This logarithm is called the principal logarithm of $A$, and denoted by $\log(A)$.
Throughout this paper, we assume that all eigenvalues of $A$ lie in the set $\mathbb{C}\setminus (-\infty,0]$, and we consider the principal logarithm of $A$.

The matrix logarithm is utilized in many fields of research, such as quantum
mechanics \cite{cosmas2007classical}, quantum chemistry
\cite{hesselmann2010random}, biomolecular dynamics
\cite{horenko2008likelihood}, buckling simulation
\cite{schenk2009modeling}, and machine learning \cite{han2015large,huang2017riemannian,don2017scalable,fitzsimons2017entropic}.
The computational methods include the inverse scaling and squaring (ISS) algorithm \cite{al2012improved}, an algorithm based on the arithmetic-geometric mean (AGM) iteration \cite{cardoso2016matrix}, and numerical quadrature.
In this paper, we focus on numerical quadrature, which employs the following integral representation (see e.g. \cite[Thm. 11.1]{higham2008functions}):
\begin{align}\label{eq:1-01}
  \log(A) = (A-I) \int_0^1 \bigl[t(A-I) + I\bigr]^{-1} \,\mathrm{d} t.
\end{align}
First, numerical quadrature can make use of the sparseness of $A$ if $A$ is sparse.
  It is useful when computing the multiplication of the matrix logarithm and a vector $\log(A)\bm{b}~(\bm{b}\in \mathbb{R}^n)$, which appears in applications such as \cite{han2015large,don2017scalable,fitzsimons2017entropic}\footnote{
    In the three applications, the computation of $\log(\det A)$ is performed based on some connections between $\log(A)\bm{b}$ and $\log(\det A)$. For more details, see, e.g., \cite{hutchinson1990stochastic}.
  }, without computing and storing dense matrices.
Conversely, the ISS algorithm and the algorithm based on the AGM iteration include the computation of the matrix square root, which means that the calculation involves dense matrices even if $A$ is sparse.
The second reason is that numerical quadrature is potentially more favorable for parallel computers because of independent computation of the integrand on each abscissa.

Because the integrand in \eqref{eq:1-01} includes matrix inversion, the computational cost of numerical quadrature depends on the number of evaluations of the integrand.
Although numerical quadrature is suitable for parallelization, the quadrature formula should be selected carefully to reduce the computational cost and save computational resources.

The method conventionally used to compute \eqref{eq:1-01} is the Gauss--Legendre (GL) quadrature.
If the spectral radius of $A-I$ is smaller than 1, the GL quadrature, which can be regarded as a rational approximation of $\log(A)$, coincides with the Pad\'{e} approximants of $\log(A)$ at $I$ \cite[Thm. 4.3]{dieci1996computational}.
Therefore, it is natural to use the GL quadrature to reduce the number of abscissas when $A$ is close to $I$.
However, the convergence of the GL quadrature becomes slow when $A$ is not close to $I$.
For example, the convergence in our experiments became slow when $A$ was ill-conditioned which may be explained by rapid changes in the integrand value when it is closer to the endpoint of the interval.

In this paper, we consider the double exponential (DE) formula \cite{takahashi1974double}, which can be used to compute integrals with singularities at one (or both) of the endpoints.
For this reason, the DE formula may be useful in scenarios in which the GL quadrature does not perform well.
However, when using the DE formula, a finite interval needs to be selected because the integrand in \eqref{eq:1-01} is transformed into a corresponding function on the infinite interval.
This selection is important, i.e., if the finite interval is too narrow, the accuracy of the computational result becomes low, but if it is too wide, the convergence of the DE formula becomes slow.

By performing an error analysis, we provide a method of selecting the appropriate finite interval, as well as two algorithms for the computation of $\log(A)$ based on the $m$-point DE formula.

The remainder of this paper is organized as follows:
in Section 2, we present an error analysis and propose two algorithms;
in Section 3, we show the results of numerical experiments;
in Section 4, we conclude the study.

\textbf{Notation:}
Unless otherwise stated, $\|\cdot\|$ denotes a consistent matrix norm, e.g., the $p$-norm or the Frobenius norm, and $\|\cdot\|_2$ and $\|\cdot\|_F$ denote the 2-norm and the Frobenius norm, respectively.
The inverse functions of $\sinh$ and $\tanh$ are referred to as $\arsinh$ and $\artanh$, respectively.

\section{Algorithms for the computation of $\log(A)$ based on the DE formula}
In this section we propose a method of selecting a finite interval for the DE formula by estimating the interval truncation error and present two algorithms.
Before considering the truncation error, let us apply the following transformations to \eqref{eq:1-01}.
By substituting $u=2t-1$ in \eqref{eq:1-01}, we obtain
\begin{align}\label{eq:2-01}
  \log(A) = (A-I)\int_{-1}^1 \left[(1+u)(A-I)+2I\right]^{-1} \,\mathrm{d} u.
\end{align}
Then, by applying the DE transformation $u = \tanh(\sinh(x))$, it follows that
\begin{align}\label{eq:2-02}
  \log(A) = (A-I)\int_{-\infty}^\infty F_{\mathrm{DE}}(x) \,\mathrm{d} x,
\end{align}
where
\begin{align}
  F_{\mathrm{DE}}(x) := \cosh(x)\sech^2(\sinh(x))\left[
      (1+\tanh(\sinh(x)))(A-I)+2I
    \right]^{-1}.
\end{align}
The matrix $(1+\tanh(\sinh(x)))(A-I)+2I$ in the integrand $F_{\mathrm{DE}}$ is nonsingular for any $x\in$ $\mathbb{R}$.
In Subsection 2.1, we derive an upper bound on the error between the integral in \eqref{eq:2-02} and the same integral defined in the finite interval $[l,r]$,
\begin{align}\label{eq:2-51}
  \left\lVert \log(A) - (A-I)\int_l^r F_{\mathrm{DE}}(x) \,\mathrm{d} x \right\rVert.
\end{align}
In Subsection 2.2, we propose a method of selecting the interval $[l,r]$ so that the relative truncation error is small than or approximately equal to a given tolerance $\epsilon$.
Our algorithms are described in Subsection 2.3.

\subsection{Estimation of the error from the interval truncation}
The error that stems from the interval truncation \eqref{eq:2-51} can be rewritten as
\begin{align}\label{eq:2-03}
  \left\lVert \log(A) - (A-I)\int_l^r F_{\mathrm{DE}}(x) \,\mathrm{d} x \right\rVert
 =  \left\lVert (A-I)\left[
   \int_{-\infty}^lF_{\mathrm{DE}}(x)\,\mathrm{d} x
   +\int_r^\infty F_{\mathrm{DE}}(x)\,\mathrm{d} x
  \right]
  \right\rVert.
\end{align}
Using the triangle inequality in the right-hand side (RHS) of \eqref{eq:2-03}, it holds that
\begin{align}\label{eq:2-04}
  & \left\lVert (A-I)\left[
    \int_{-\infty}^lF_{\mathrm{DE}}(x)\,\mathrm{d} x
    +\int_r^\infty F_{\mathrm{DE}}(x)\,\mathrm{d} x
   \right]
   \right\rVert \\
  &\quad \le  \left\lVert (A-I)\int_{-\infty}^lF_{\mathrm{DE}}(x)\,\mathrm{d} x \right\rVert
  +  \left\lVert (A-I)\int_r^\infty F_{\mathrm{DE}}(x)\,\mathrm{d} x \right\rVert .
\end{align}
By estimating the RHS of \eqref{eq:2-04}, we obtain an upper bound on \eqref{eq:2-51}.

Initially, we focus on the first term which is on the RHS of \eqref{eq:2-04}.
To avoid cumbersome notation, instead of $F_{\mathrm{DE}}$, which including hyperbolic functions, we consider the integrand of \eqref{eq:1-01}.
By applying the transformation $t=[\tanh(\sinh(x))+1]/2$, we have
\begin{align} \label{eq:2-05}
  \left\lVert (A-I)\int_{-\infty}^lF_{\mathrm{DE}}(x)\,\mathrm{d} x \right\rVert
  = \left\lVert (A-I)\int_0^a F(t) \,\mathrm{d} t \right\rVert,
\end{align}
where,
\begin{align}
  F(t):=[t(A-I)+I]^{-1},\quad a:=\frac{\tanh(\sinh(l))+1}{2}.
\end{align}
The following lemma shows an upper bound on the RHS of \eqref{eq:2-05}, if $a$ is small enough to warrant the use of the Neumann series expansion of $F(t)$.
\begin{lem}\label{thrm:error_l}
  Suppose that $A\ne I$.
  Then, for $a\in \bigl(0, 1/(2\|A-I\|)\bigr]$, we have
  \begin{align}\label{eq:2-06}
    \left\lVert (A-I)\int_0^a F(t) \,\mathrm{d} t \right\rVert  \le \frac{3a}{2} \|A-I\|.
  \end{align}
\end{lem}
\begin{pf}
  For all $t\in[0,a]$ where $a\in \bigl(0, 1/(2\|A-I\|)\bigr]$, it holds that
  \begin{align}
    \left\lVert t(I-A) \right\rVert    \le \frac{1}{2}~ (<1).
  \end{align}
  By applying the Neumann series expansion to $F(t)$ we get the following:
  \begin{align}
    F(t) = [t(A-I)+I]^{-1} = \sum_{k=0}^\infty t^k(I-A)^k.
  \end{align}
  Therefore, the integral of $F$ can be rewritten as follows:
  \begin{align}
    (A-I)\int_0^a F(t) \,\mathrm{d} t
    &= (A-I)\int_0^a \left[\sum_{k=0}^\infty t^k(I-A)^k\right] \,\mathrm{d} t\label{eq:2-50}\\
    &= (A-I)\sum_{k=0}^\infty\frac{a^{k+1}}{k+1}(I-A)^k\\
    &= a(A-I) + a(A-I)\sum_{k=1}^\infty\left[\frac{a^k}{k+1}(I-A)^k\right].
  \end{align}
  By using triangle inequality and consistency of the norm we get the following:
  \begin{align}
    \left\lVert (A-I)\int_0^a F(t) \,\mathrm{d} t \right\rVert
   &\le a \left\lVert A-I \right\rVert  + a \left\lVert A-I \right\rVert \sum_{k=1}^\infty\left[
     \frac{1}{k+1} \left\lVert a(A-I) \right\rVert ^k
   \right]\\
   &\le a \left\lVert A-I \right\rVert  + a \left\lVert A-I \right\rVert \sum_{k=1}^\infty\left[
     \frac{1}{2} \left(\frac{1}{2}\right)^k
   \right]\\
   &= \frac{3a}{2} \left\lVert A-I \right\rVert .
  \end{align}
  \hfill $\square$
\end{pf}

The calculation to estimate the second term on the RHS of \eqref{eq:2-04} is similar to that of the first term.
By applying the transformation $t=[\tanh(\sinh(x))+1]/2$, we get the following:
\begin{align}\label{eq:2-08}
  \left\lVert (A-I)\int_r^\infty F_{\mathrm{DE}}(x) \,\mathrm{d} x \right\rVert  =  \left\lVert (A-I)\int_b^1 F(t) \,\mathrm{d} t \right\rVert ,
\end{align}
where $b=[\tanh(\sinh(r))+1]/2$.

The following lemma shows an upper bound on \eqref{eq:2-08}, if $b$ is close enough to 1 to warrant the use of the Neumann series expansion of $F(t)$.
\begin{lem}\label{thrm:error_r}
  For $b \in \bigl[2 \|A^{-1}\|/(2\|A^{-1}\|+1), 1\bigr)$, we have
  \begin{align}\label{eq:2-10}
    \left\lVert (A-I)\int_b^1 F(t) \,\mathrm{d} t \right\rVert  \le \left(-\log(b)+\frac{1-b}{2b}\right) \|A-I\|\|A^{-1}\|.
    \end{align}
\end{lem}
\begin{pf}
  The outline of this proof is similar to that of Lemma \ref{thrm:error_l}.
  For all $t\in[b,1]$ where $b \in \bigl[2 \|A^{-1}\|/(2\|A^{-1}\|+1), 1\bigr)$, it holds that
  \begin{align}
    \left\lVert \frac{t-1}{t}A^{-1} \right\rVert  \le \frac{1}{2}~(<1).
  \end{align}
  By applying the Neumann series expansion to $F(t)$, we get
  \begin{align}
    F(t) &= [t(A-I)+I]^{-1}\\
    &= \frac{1}{t}A^{-1}\sum_{k=0}^\infty\left(\frac{t-1}{t}\right)^kA^{-k}.
  \end{align}
  Therefore, the integral of $F$ can be rewritten as follows:
  \begin{align}
    (A-I)\int_b^1 F(t) \,\mathrm{d} t
    &= (A-I)\int_b^1
      \left[\frac{1}{t}A^{-1}\sum_{k=0}^\infty\left(\frac{t-1}{t}\right)^kA^{-k}\right]
     \,\mathrm{d} t\\
    &= (A-I)A^{-1}\left[
      -\log(b)I+\sum_{k=1}^\infty A^{-k}\int_b^1 \frac{1}{t} \left(\frac{t-1}{t}\right)^k \,\mathrm{d} t
    \right].
  \end{align}
  By using the triangle inequality and consistency of the norm, it follows
  \begin{align} \label{eq:2-07}
    \left\lVert (A-I)\int_b^1 F(t) \,\mathrm{d} t \right\rVert
     &\le  \left\lVert A-I \right\rVert  \left\lVert A^{-1} \right\rVert \left\{
       |-\log(b)| + \sum_{k=1}^\infty\left[
          \left\lVert A^{-1} \right\rVert ^k \int_b^1 \frac{1}{t} \left\lvert\frac{t-1}{t}\right\rvert^k \,\mathrm{d} t
       \right]
     \right\}.
  \end{align}
  In \eqref{eq:2-07}, it holds that $|-\log(b)|=-\log(b)$ because $b\in\bigl[2 \|A^{-1}\|/(2\|A^{-1}\|+1), 1\bigr)$, and $|(t-1)/t| = (1-t)/t$ because $t\in[b,1]$.
  For the second term in the bracket on the RHS of \eqref{eq:2-07}, we have:
  \begin{align}
    \sum_{k=1}^\infty\left[ \left\lVert A^{-1} \right\rVert^k \int_b^1 \frac{1}{t} \left(\frac{1-t}{t}\right)^k \,\mathrm{d} t\right]
    &\le\sum_{k=1}^\infty\left[ \left\lVert A^{-1} \right\rVert^k \int_b^1 \frac{1}{b} \left(\frac{1-t}{b}\right)^k \,\mathrm{d} t\right]\\
    &= \frac{1-b}{b}\sum_{k=1}^\infty\left[\frac{1}{k+1} \left\lVert \frac{1-b}{b}A^{-1} \right\rVert ^k\right]\\
    &\le \frac{1-b}{b} \sum_{k=1}^\infty\left[\frac{1}{2}\left(\frac{1}{2}\right)^k\right]\\
    & = \frac{1-b}{2b}. \label{eq:2-09}
  \end{align}
  By substituting \eqref{eq:2-09} in \eqref{eq:2-07}, we obtain the inequality \eqref{eq:2-10}.
  \hfill $\square$
\end{pf}

In the final part of this Subsection, we estimate the upper bound on \eqref{eq:2-51}.
\begin{prop}\label{thrm:2-01}
  Suppose that $A\ne I$.
  For a given interval $[l,r]$, let $a=[\tanh(\sinh(l))+1]/2$ and $b=[\tanh(\sinh(r))+1]/2$.
  Then, if $a \le 1/(2\|A-I\|)$ and $b \ge 2\|A^{-1}\|/(2\|A^{-1}\|+1)$, it holds that
  \begin{align}\label{eq:2-111}
    \left\lVert \log(A) - (A-I)\int_l^r F_{\mathrm{DE}}(x) \,\mathrm{d} x \right\rVert
    \le \frac{3}{2}a\|A-I\| + \left(-\log(b)+\frac{1-b}{2b}\right)\|A-I\|\|A^{-1}\|.
  \end{align}
\end{prop}
\begin{pf}
  By combining \eqref{eq:2-03} and \eqref{eq:2-04}, and
  as well as substituting \eqref{eq:2-06} and \eqref{eq:2-10} in \eqref{eq:2-04}, we get \eqref{eq:2-111}.
  \hfill $\square$
\end{pf}

\subsection{Setting the integration interval}
To develop algorithms for computing $\log(A)$ based on the DE formula, we need to determine the appropriate finite integration interval, $[l,r]$ in advance.
The finite interval should be ideally set so that the relative error is guaranteed to be smaller than or equal to a  given tolerance, $\epsilon >0$, i.e.,
\begin{align}\label{eq:2-131}
  \frac{
  \left\lVert \log(A) - (A-I)\int_l^r F_{\mathrm{DE}}(x) \,\mathrm{d} x \right\rVert
  }{
  \|\log(A)\|
  } \le  \epsilon.
\end{align}
To accomplish this, a lower bound on $\|\log (A)\|$ must be estimated.
The following lemma shows a lower bound in terms of the spectral radius of $A$.

\begin{lem}\label{thrm:log_lower_bound}
  Let $\rho(\cdot)$ be the spectral radius, i.e., the largest absolute value of eigenvalues.
  Then, the following two inequalities hold:
  \begin{align}
    &\|\log(A)\| \ge |\log(\rho(A))|, \label{eq:2-21}\\
    &\|\log(A)\| \ge |\log(\rho(A^{-1}))|. \label{eq:2-22}
  \end{align}
\end{lem}
\begin{pf}
  By using the consistency of the norm and the fact any eigenvalue of $\log(A)$ is equal to $\log(\lambda)$, for some $\lambda$ being an eigenvalue of $A$, we have the lower bound in \eqref{eq:2-21} as $\|\log(A)\|\ge\rho(\log(A))\ge|\log(\rho(A))|$.
  The lower bound in \eqref{eq:2-22} can be obtained in a similar way.
  \hfill $\square$
\end{pf}

In the following proposition, we show how to set a finite interval such that the relative truncation error in 2-norm is smaller than or equal to a  given tolerance, $\epsilon>0$.
\begin{prop}\label{thrm:settings_of_lr}
  Suppose that $A\ne I$.
  Let $\theta$ be a lower bound on $\|\log(A)\|_2$, the tolerance $\epsilon>0$ satisfy
  \begin{align}\label{eq:2-151}
    \epsilon < \frac{
      3\|A-I\|_{2}\|A^{-1}\|_{2}
    }{
      \theta(1+\|A^{-1}\|_{2})
    },
  \end{align}
  and $s$ be a real positive solution of the equation
  \begin{align}\label{eq:2-133}
    \frac{1}{\theta}\left(-\log(s)+\frac{1-s}{2s}\right)\|A-I\|_2\|A^{-1}\|_2
    = \frac{\epsilon}{2}.
  \end{align}
  If
  \begin{align}
    l := \arsinh(\artanh(2a - 1)),\quad r := \arsinh(\artanh(2b - 1)),
  \end{align}
  where
  \begin{align}\label{eq:2-132}
    a:= \min\left\{\frac{\theta\epsilon}{3\|A-I\|_2},\frac{1}{2\|A-I\|_2}\right\},
    \quad
    b:=\max\left\{s,\frac{2\|A^{-1}\|_2}{2\|A^{-1}\|_2+1}\right\},
  \end{align}
  then, it holds that
  \begin{align}\label{eq:2-30}
    \frac{
      \left\lVert
        \log(A) - (A-I)\int_l^r F_{\mathrm{DE}}(x) \,\mathrm{d} x
      \right\rVert_2
    }{
      \|\log(A)\|_2
    } \le \epsilon.
  \end{align}
\end{prop}
\begin{pf}
  From the definition of $a$ and $b$, it is true that $a\le 1/(2\|A-I\|_2)$ and $b\ge 2\|A^{-1}\|_2/(2\|A^{-1}\|_2+1)$.
  In addition, from
  \begin{align}
    b-a
    &\ge \frac{2\|A^{-1}\|_2}{2\|A^{-1}\|_2+1} - \frac{\theta\epsilon}{3\|A-I\|_2}
    > \frac{2\|A^{-1}\|_2}{2\|A^{-1}\|_2+1} - \frac{\|A^{-1}\|_2}{1+\|A^{-1}\|_2}\\
    &=\frac{2\|A^{-1}\|_2(1+\|A^{-1}\|_2) -(2\|A^{-1}\|_2+1)\|A^{-1}\|_2 }{(2\|A^{-1}\|_2+1)(1+\|A^{-1}\|_2)}\\
    &=\frac{\|A^{-1}\|_2}{(2\|A^{-1}\|_2+1)(1+\|A^{-1}\|_2)}
    > 0,
  \end{align}
  it follows that $a<b$, and $l<r$.
  Therefore, we can choose a finite interval $[l,r]$ that satisfies the assumptions of Proposition \ref{thrm:2-01}.
  By dividing the inequality \eqref{eq:2-111} by $\|\log(A)\|_2\ge \theta$, it follows:
  \begin{align}\label{eq:2-31}
    \frac{ \left\lVert \log(A) - (A-I)\int_l^r F_{\mathrm{DE}}(x) \,\mathrm{d} x \right\rVert _2}{\|\log(A)\|_2}
    \le \frac{3a}{2\theta}\|A-I\|_2+\frac{1}{\theta}\left(-\log(b)+\frac{1-b}{2b}\right)\|A-I\|_2\|A^{-1}\|_2.
  \end{align}
  From the definition of $a$ and $b$, it holds that $a\le \theta\epsilon/(3\|A-I\|_2)$ and $b\ge s$.
  Therefore,
  \begin{align}\label{eq:2-32}
    \frac{3a}{2\theta}\|A-I\|_2 \le \frac{\epsilon}{2},\quad
    \frac{1}{\theta}\left(-\log(b)+\frac{1-b}{2b}\right)\|A-I\|_2\|A^{-1}\|_2 \le \frac{\epsilon}{2}.
  \end{align}
  We obtain \eqref{eq:2-30} by substituting \eqref{eq:2-32} in \eqref{eq:2-31}.
  \hfill $\square$
\end{pf}

\subsection{Algorithms}
In this subsection we establish two algorithms based on the results in subsections 2.1 and 2.2.
One of the algorithms is designed to compute $\log(A)$ using the $m$-point DE formula on a finite interval with an interval truncation error smaller than or approximately equal to a  given tolerance, $\epsilon>0$.
The other algorithm is an adaptive quadrature algorithm designed to compute $\log(A)$ by automatically adding abscissas until the error is smaller than or approximately equal to a given tolerance $\zeta>0$.

If the tolerance $\epsilon$ given in Proposition \ref{thrm:settings_of_lr} is sufficiently small, a linear approximation to the nonlinear equation \eqref{eq:2-133} can be used to determine an appropriate interval. We describe our calculation in detail below.

Suppose that $\epsilon$ is sufficiently small and the solution $s$ of \eqref{eq:2-133} is approximately equal to 1.
Then, because $3(1-s)/2$ is the first-order Taylor approximation to $-\log(s)+(1-s)/2s$ at $s =1$, the solution $s$ can be approximated by using the solution $\tilde{s}$ of the following equation:
\begin{align}\label{eq:2-140}
  \frac{1}{\theta} \frac{3(1-\tilde{s})}{2} \|A-I\|_2\|A^{-1}\|_2
  = \frac{\epsilon}{2}.
\end{align}
The solution of \eqref{eq:2-140} is given by
\begin{align}
  \tilde{s} =  1-\frac{\theta\epsilon}{3\|A-I\|_2\|A^{-1}\|_2}.
\end{align}

Under the assumptions of Proposition \ref{thrm:settings_of_lr}, by choosing $\tilde{b}$ as
\begin{align}
  \tilde{b} = \max \left\{
    \tilde{s},~ \frac{2\|A^{-1}\|_2}{2\|A^{-1}\|_2+1}
    \right\}
\end{align}
instead of \eqref{eq:2-132}, and setting $\tilde{r}=\arsinh(\artanh(2\tilde{b}-1))$, the interval truncation will be smaller than or approximately equal to $\epsilon$.
The summary of the first algorithm, which computes $\log(A)$ using the $m$-point DE formula whose interval truncation error is smaller than or approximately equal to $\epsilon$ is as shown in Algorithm \ref{alg:1}.

\begin{algorithm}[htbp]
  \caption{Computation of $\log(A)$ based on the DE formula.}
  \begin{algorithmic}[1]
    \Statex \textbf{Input:} $A\in\mathbb{R}^{n\times n}$, $m\in\mathbb{N}$, $\epsilon>0$ a tolerance for the interval truncation error
    \Statex \textbf{Output:} $X\approx \log(A)$
    \State Set $F_{\mathrm{DE}}(x)=\cosh(x)\sech^2(\sinh(x))\left[
      (1+\tanh(\sinh(x)))(A-I)+2I
    \right]^{-1}$
    \State Compute $\|A-I\|_2$, $\|A^{-1}\|_2$, $\rho(A)$
    \State $\theta = |\log(\rho(A))|$
    \State $\displaystyle\epsilon_{\mathrm{max}} = \frac{3}{\theta} \frac{\|A-I\|_2\|A^{-1}\|_2}{1+\|A^{-1}\|_2}$
    \If {$\displaystyle \epsilon \ge \epsilon_{\mathrm{max}}$}
      \State $\displaystyle \epsilon \gets \epsilon_{\mathrm{max}}/2$
    \EndIf
    \State $\displaystyle a = \min \left\{
      \frac{\theta\epsilon}{3\|A-I\|_2}, \frac{1}{2\|A-I\|_2}
    \right\}$
    \State $\displaystyle b = \max \left\{
      1-\frac{\theta\epsilon}{3\|A-I\|_2\|A^{-1}\|_2}, \frac{2\|A^{-1}\|_2}{2\|A^{-1}\|_2+1}
    \right\}$
    \State $l = \arsinh(\artanh(2a - 1))$
    \State $r = \arsinh(\artanh(2b - 1))$
    \State $h = (r-l)/(m-1)$
    \State $\displaystyle T = \frac{h}{2}(F_{\mathrm{DE}}(l)+F_{\mathrm{DE}}(r)) + h\sum_{i=1}^{m-2}F_{\mathrm{DE}}(l+ih)$
    \State $X = (A-I)T$
  \end{algorithmic}
  \label{alg:1}
\end{algorithm}

When $\epsilon$ is sufficiently small, an accurate computation of $\|I-A\|_2$, $\|A^{-1}\|_2$ and $\rho(A)$ at Step 2 of Algorithm \ref{alg:1} may not be required because the errors that stem from of these values have little effect on the accuracy of $\log(A)$. We give more detail in the following paragraph.

Assume that $\epsilon$ in Algorithm \ref{alg:1} is sufficiently small, and that $a$ and $b$ in Step 9 are chosen as
\begin{align}
  a = \frac{\theta\epsilon}{3\|A-I\|_2},\quad
  b = 1-\frac{\theta\epsilon}{3\|A-I\|_2\|A^{-1}\|_2},
\end{align}
where $\theta$ is a lower bound on $\|\log(A)\|_2$.
Let $\Delta_1$ and $\Delta_2$ be the relative errors of $\theta/\|A-I\|_2$ and $1/\|A^{-1}\|_2$ respectively.
Then, the computational result of $a$ is equal to
\begin{align}
  \frac{\epsilon}{3}\frac{\theta}{\|A-I\|_2}(1+\Delta_1)
  = \frac{\theta}{3\|A-I\|_2} \epsilon(1+\Delta_1),
\end{align}
and that of $b$ is equal to
\begin{align}
  1-\frac{\epsilon}{3}\frac{\theta}{\|A-I\|_2}(1+\Delta_1)
  \frac{1}{\|A^{-1}\|_2}(1+\Delta_2)
  = 1-\frac{\theta}{3\|A-I\|_2\|A^{-1}\|_2}\epsilon(1+\Delta_1+\Delta_2+\Delta_1\Delta_2).
\end{align}
Therefore, the upper bound on the truncation error $\epsilon$ computed by considering the relative errors $\Delta_1,~\Delta_2$ is almost equal to the upper bound on the truncation error when the tolerance is set as $\epsilon(1+\Delta_1+\Delta_2+\Delta_1\Delta_2)$.
For example, when $\Delta_1,\Delta_2=10^{-2}$, $\epsilon(1+\Delta_1+\Delta_2+\Delta_1\Delta_2)\approx 1.02\epsilon$, which means that the upper bound on the truncation error changes by approximately 2\%.
If $\epsilon$ is sufficiently small, the effect of these errors will be negligible.

At Step 3, a lower bound on $\|\log(A)\|_2$ is computed based on \eqref{eq:2-21}.
By setting $\theta=\max\{|\log(\rho(A))|,$ $|\log(\rho(A^{-1}))|\}$, a tighter lower bound can be obtained.
In particular, when $A$ is positive definite, $|\log(\rho(A^{-1}))|$ can be obtained without additional computation because $\rho(A^{-1})=\|A^{-1}\|_2$ is already computed in Step 2.

The computational cost of Algorithm 1 for dense $A$ is $(2m + 2) n^3 + \mathcal{O}(n^2)$.
When $A$ is sparse, evaluating $\log(A) \bm{b}$ using Algorithm 1 has computational cost $mc_{\mathrm{abscissa}} + c_{\mathrm{mul}} + c_{\mathrm{param}}$, where $c_{\mathrm{abscissa}}$ is the computational cost of computing $F_{\mathrm{DE}}(x)\bm{b}$,
$c_{\mathrm{mul}}$ is the computational cost of a matrix-vector multiplication, and $c_{\mathrm{param}}$ is the computational cost of computing parameters $\rho(A)$, $\|A-I\|_2$, and $\|A^{-1}\|_2$.
If the parameters are computed approximately and $F_{\mathrm{DE}}(x)\bm{b}$ is computed accurately, then $c_{\mathrm{param}}$ will be smaller than $c_{\mathrm{abscissa}}$.
Therefore, the computational cost will largely depend on $m$.

Once an appropriate finite interval is obtained, the accuracy of the DE formula can be improved with the following procedure.
Let $m$ be the number of abscissas, $h=(r-l)/(m-1)$ be the mesh size, and $T(h)$ be the trapezoidal rule for the mesh size $h$:
\begin{align}
  T(h) := \frac{h}{2} \left(F_{\mathrm{DE}}(l) + F_{\mathrm{DE}}(r)\right) + h \sum_{i=1}^{m-2} F_{\mathrm{DE}}(l+ih).
\end{align}
Then, $T(h/2)$ can be computed using $T(h)$:
\begin{align}
  T\left(\frac{h}{2}\right) = \frac{1}{2}T(h) + \frac{h}{2}\sum_{i=1}^{m-1} F_{\mathrm{DE}}\left(l+(2i-1)\frac{h}{2}\right).
\end{align}
In addition, we can apply the following estimation of the trapezoidal error for a sufficiently small value of $h$ using \cite[Eq. (4.3)]{cardoso2012computation}:
\begin{align}\label{eq:2-40}
  \left\lVert \int_l^r F_{\mathrm{DE}}(x) - T\left(\frac{h}{2}\right) \right\rVert
  \approx \frac{1}{3} \left\lVert T(h) - T\left(\frac{h}{2}\right) \right\rVert.
\end{align}
Our adaptive quadrature algorithm which is based on \eqref{eq:2-40} is presented as Algorithm \ref{alg:2}.
\begin{algorithm}
  \caption{Computation of $\log(A)$ by adaptive quadrature based on the DE formula.}
  \begin{algorithmic}[1]
    \Statex \textbf{Input:} $A\in\mathbb{R}^{n\times n}$, $m_0\in\mathbb{N}$, $\epsilon>0$ a tolerance for the interval truncation error, $\zeta>0$ a tolerance for the trapezoidal truncation error.
    \Statex \textbf{Output:} $X\approx \log(A)$
    \State Set $F_{\mathrm{DE}}(x)=\cosh(x)\sech^2(\sinh(x))\left[
      (1+\tanh(\sinh(x)))(A-I)+2I
    \right]^{-1}$
    \State Computing $l,~r,~\theta$ using steps 2 to 11 of Algorithm 1
    \State $h_0 = (r-l) / (m_0-1)$
    \State $T_0 = \frac{h_0}{2}F_{\mathrm{DE}}(l) + \frac{h_0}{2}F_{\mathrm{DE}}(r)
      + h_0\sum_{i=1}^{m_0-2} F_{\mathrm{DE}}(l+ih)$
    \For{$k=0,1,2,\ldots$ until convergence}
      \State $h_{k+1}=h_k/2$
      \State $T_{k+1}=\frac{1}{2}T_k + h_{k+1}\sum_{i=1}^{m_k-1} F_{\mathrm{DE}}(l+(2i-1)h_{k+1})$
      \State $m_{k+1}=2m_k-1$
      \If {$\frac{1}{3}\|T_{k+1}-T_k\|/\theta \le \zeta$}
        \State $T=T_{k+1}$
        \State \textbf{break}
      \EndIf
    \EndFor
    \State $X = (A-I)T$
  \end{algorithmic}
  \label{alg:2}
\end{algorithm}

The computational cost of Algorithm \ref{alg:2} for dense $A$ is $(2m_{k+1}+2)n^3 + \mathcal{O}(n^2)=[2^{k+1}(m_0-1)+4]n^3 + \mathcal{O}(n^2)$, and the computational cost of $\log(A)\bm{b}$ with sparse $A$ is $[2^k(m_0-1)+1]c_{\mathrm{abscissa}} + c_{\mathrm{mul}} + c_{\mathrm{param}}$.

\section{Numerical experiments}
The numerical experiments were carried out using Julia 1.0 on a Core-i7 (3.4GHz) CPU with 16GB RAM.
We used the IEEE double precision arithmetic.
We computed abscissas and weights in the GL quadrature with \texttt{QuadGK.jl} (https://github.com/JuliaMath/QuadGK.jl).

\subsection{Experiment 1: Checking the convergence}
In this experiment, we checked the convergence of the GL quadrature and the DE formula.
Our test matrices are presented in Table \ref{tab:test_matrices}.
\begin{table}[htbp]
  \caption{Test matrices. The condition number of $A$ is denoted by $\kappa_2(A)=\|A\|_2\|A^{-1}\|_2$.}
  \centering
  \begin{tabular}{lrll}
    \hline
    Matrix & \multicolumn{1}{l}{Size} & $\kappa_2(A)$ & Structure \\ \hline
    \texttt{SPD 1} & 50 & $1.0\times 10^1$ & SPD\\
    \texttt{SPD 2} & 50 & $1.0\times 10^4$ & SPD\\
    \texttt{SPD 3} & 50 & $1.0\times 10^7$ & SPD\\
    \texttt{parter} \cite{zhang2016matrix} & 10 & $2.4\times 10^{0}$ & Nonsymmetric \\
    \texttt{frank} \cite{zhang2016matrix} & 10 & $2.9\times 10^{7}$ & Nonsymmetric \\
    \texttt{vand} \cite{zhang2016matrix} & 10 & $3.1\times 10^{12}$ & Nonsymmetric \\
    \texttt{bcsstk02} \cite{davis2011university} & 66 & $4.3\times 10^{3}$ & SPD \\
    \texttt{bcsstk03} \cite{davis2011university} & 112 & $6.8\times 10^{6}$ & SPD \\
    \texttt{ck104} \cite{davis2011university} & 104 & $5.5\times 10^{3}$ & Nonsymmetric \\ \hline
  \end{tabular}
  \label{tab:test_matrices}
\end{table}
We generated the first three matrices in Table \ref{tab:test_matrices} by using the following procedure:
\begin{enumerate}
  \item We generated an orthogonal matrix $Q$ by QR decomposition of a random $50\times 50$ matrix.
  \item We generated a diagonal matrix whose diagonal elements were from the geometric sequence: $\{d_i\}_{i=1,\ldots,50}$ where $d_1=\kappa^{-1/2}$ and $d_{50}=\kappa^{1/2}$ for $\kappa=10^1$.
  \item $A = Q D Q^\top$.
  \item We repeated Step 2 and Step 3 by setting $\kappa$ as $10^4$ and $10^7$.
\end{enumerate}

The experimental procedure is as follows:
\begin{enumerate}
  \item We scaled the test matrices as $\tilde{A} = (10/\rho(A)) A$ because some matrices had values that were too large to use in computation.
  \item We computed the reference $\log(\tilde{A})$ with the arbitrary precision arithmetic and rounded the result to double precision. We implemented the ISS algorithm \cite[Alg. 11.10]{higham2008functions} with the \texttt{BigFloat} type of Julia.
  \item We computed $\log(\tilde{A})$ using Algorithm \ref{alg:1}, where $\epsilon=2^{-53}\approx 1.1\times 10^{-16}$. If the test matrix was symmetric positive definite, we set $\theta=\max\{|\log(\rho(\tilde{A}))|,|\log(\rho(\tilde{A}^{-1}))|\}$ as stated in Subsection 2.3.
  We computed $\rho(\tilde{A})$ using the \texttt{eigvals} function of Julia, which computes all eigenvalues of $\tilde{A}$.
  Similarly, we computed $\|I-\tilde{A}\|_2$ and $\|\tilde{A}^{-1}\|_2$ using the \texttt{svdvals} function, which computes all singular values of $\tilde{A}$
  \footnote{
    If $A$ is large, using \texttt{eigvals} and \texttt{svdvals} may be inefficient.
    Instead, for Julia, a package \texttt{Arapack.jl} (https://github.com/JuliaLinearAlgebra/Arpack.jl) is available, which can compute a part of eigenvalues and singular values based on the Lanczos (or the Arnoldi) process with the desired accuracy.
    We present some numerical results, for which $\rho(\tilde{A})$, $\|\tilde{A}-I\|_2$, and $\|\tilde{A}^{-1}\|_2$ are computed with low accuracy, as shown in Appendix A.
  }.
  \item We computed $\log(\tilde{A})$ by applying the GL quadrature to \eqref{eq:2-01}.
\end{enumerate}

\begin{figure}[htbp]\centering
  \includegraphics[width=0.95\linewidth]{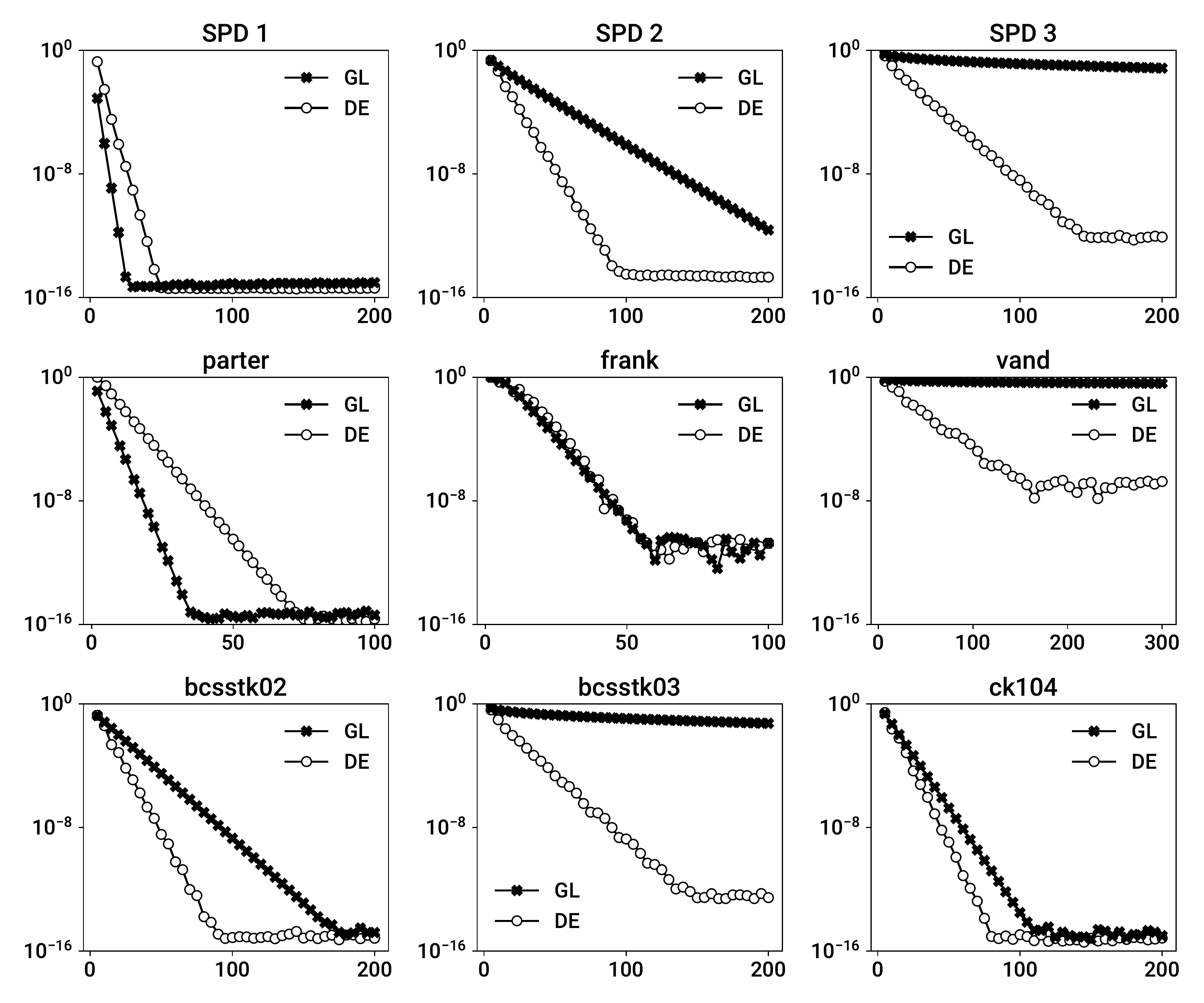}
  \caption{Convergence histories of the DE formula (Algorithm \ref{alg:1}) and the GL quadrature. The vertical axes show the relative error, $\|\log(\tilde{A})-X\|_F/\|\log(\tilde{A})\|_F$, and the horizontal axes show the number of abscissas, $m$.}
  \label{fig:1}
\end{figure}
Figure \ref{fig:1} shows the convergence histories of the DE formula and the GL quadrature for each matrix.
Several observations can be made:
\begin{itemize}
  \item The accuracy of the DE formula is almost the same as that of the GL quadrature, and the accuracies of the DE formula and the GL quadrature depend on the condition number of test matrices.
  \item For well-conditioned matrices, such as \texttt{SPD 1} and \texttt{parter} matrix, the GL quadrature converged faster than the DE formula.
  Conversely, for the ill-conditioned matrices, such as \texttt{SPD 3} and \texttt{vand} matrix, the DE formula converged faster than the GL quadrature.
\end{itemize}

The above observations suggest that Algorithm \ref{alg:1} selects an appropriate interval and the DE formula is suitable for ill-conditioned matrices.

\subsection{Experiment 2: Checking Algorithm 2}
In this experiment, we check the performance of Algorithm \ref{alg:2} by using the same matrices that were used in Experiment 1 (see Subsection 3.1).
We compared Algorithm \ref{alg:2} with Algorithm \ref{alg:3}, which is based on the GL quadrature (see Appendix B).

We conducted the experiment using the following procedure:
\begin{enumerate}
  \item We computed $\log(\tilde{A})$ by using Algorithm \ref{alg:2}.
  We set $m_0=16$ and $\zeta=\epsilon\in\{10^{-8},10^{-11}\}$.
  In Step 2 of Algorithm \ref{alg:2}, which calls Algorithm 1, we computed the spectral radius and the 2-norm of matrices using \texttt{eigvals} and \texttt{svdvals} functions, as is done in Experiment 1.
  We stopped the computation once the number of integrand evaluations reached 1921.
  \item We computed $\log(\tilde{A})$ by using Algorithm \ref{alg:3}.
  We set $m_0=16$ and $\zeta\in\{10^{-8},10^{-11}\}$.
  In Step 2 of Algorithm \ref{alg:3}, the lower bound $\theta$ was computed using \eqref{eq:2-21} and \eqref{eq:2-22} in the same way as was done for the DE formula.
  If the number of integrand evaluations was more than 2032, we stopped the computation.
\end{enumerate}

\begin{table}[htbp]
  \centering
  \caption{Comparison between Algorithm \ref{alg:2} (based on the DE formula) and Algorithm \ref{alg:3} (based on the GL quadrature), in terms of the number of integrand evaluations (in bold) and the relative error when the algorithm stopped (in parentheses). The notation ``\textbf{--} (--)'' means that the algorithms did not stop before the number of integrand evaluations reached the limit.}
  \begin{tabular}{l|ll|ll}\hline
    & \multicolumn{2}{l|}{Algorithm \ref{alg:2} (DE)} & \multicolumn{2}{l}{Algorithm \ref{alg:3} (GL)} \\
    $\zeta$&  $10^{-8}$ &  $10^{-11}$ &  $10^{-8}$ &  $10^{-11}$ \\\hline
    \texttt{SPD 1} &  \textbf{61} ($2.2 \times 10^{-9}$) &  \textbf{61} ($2.7 \times 10^{-12}$) &  \textbf{48} ($4.6 \times 10^{-16}$) &  \textbf{112} ($5.7 \times 10^{-16}$) \\
    \texttt{SPD 2} &  \textbf{121} ($6.7 \times 10^{-10}$) &  \textbf{241} ($6.4 \times 10^{-13}$) &  \textbf{1008} ($1.8 \times 10^{-15}$) &  \textbf{1008} ($1.8 \times 10^{-15}$) \\
    \texttt{SPD 3} &  \textbf{241} ($3.0 \times 10^{-10}$) &  \textbf{481} ($4.9 \times 10^{-13}$) &  \textbf{--} (--) &  \textbf{--} (--) \\
    \texttt{parter} &  \textbf{61} ($2.6 \times 10^{-9}$) &  \textbf{121} ($2.3 \times 10^{-12}$) &  \textbf{112} ($3.3 \times 10^{-16}$) &  \textbf{112} ($3.3 \times 10^{-16}$) \\
    \texttt{frank} &  \textbf{481} ($1.0 \times 10^{-12}$) &  \textbf{1921} ($2.1 \times 10^{-13}$) &  \textbf{496} ($1.5 \times 10^{-11}$) &  \textbf{--} (--) \\
    \texttt{vand}&  \textbf{--} (--) & \textbf{--} (--) &  \textbf{--} (--) &  \textbf{--} (--) \\
    \texttt{bcsstk02} &  \textbf{121} ($2.8 \times 10^{-9}$) &  \textbf{121} ($3.1 \times 10^{-12}$) &  \textbf{496} ($1.7 \times 10^{-15}$) &  \textbf{1008} ($1.0 \times 10^{-15}$) \\
    \texttt{bcsstk03} &  \textbf{241} ($1.4 \times 10^{-9}$) &  \textbf{241} ($1.5 \times 10^{-12}$) &  \textbf{--} (--) &  \textbf{--} (--) \\
    \texttt{ck104} &  \textbf{121} ($6.7 \times 10^{-10}$) &  \textbf{121} ($7.4 \times 10^{-13}$) &  \textbf{496} ($2.2 \times 10^{-15}$) &  \textbf{496} ($2.2 \times 10^{-15}$) \\ \hline
  \end{tabular}
  \label{tab:result2}
\end{table}

The number of integrand evaluations and the corresponding relative error when the two algorithms stopped are shown in Table \ref{tab:result2}.

Several observations can be made:
\begin{itemize}
  \item Algorithm \ref{alg:2} successfully computed the logarithm with the desired accuracy within 2000 integrand evaluations for all test matrices, except \texttt{vand} matrix, whereas Algorithm \ref{alg:3} could not succeed for three matrices.
  For \texttt{vand} matrix, the stopping criterion $\zeta$ may be too strict because the convergence history of \texttt{vand} matrix (as shown in Figure \ref{fig:1}) hardly reach the value of $10^{-8}$.
  Our future studies will focus on the method for selecting a suitable stopping criterion.
  \item Even if the condition number of $A$ is small as is the case for \texttt{SPD 1} and \texttt{parter} matrix, the number of integrand evaluations of Algorithm \ref{alg:2} could be smaller than that of Algorithm \ref{alg:3} because Algorithm \ref{alg:2} can reuse all previous results when improving accuracy.
\end{itemize}

These observations show that Algorithm \ref{alg:2} can be a practical choice for the computation of the matrix logarithm by numerical quadrature.

\section{Conclusion}
In this paper, we focused on the DE formula as a new choice for the numerical quadrature formula of $\log(A)$.
In order to utilize the DE formula, we proposed a method for selecting an appropriate finite interval based on error analysis, and we proposed two algorithms for practical computation.

We carried out two numerical experiments.
The first experiment showed that the DE formula converged faster than the GL quadrature for ill-conditioned matrices.
The second experiment demonstrated that the proposed adaptive quadrature algorithm worked well with appropriate stopping criteria.

Our future work will focus on three problems.
The first one is the analyses of the convergence rate for the DE formula and the GL formula, the second one is a method of selecting appropriate stopping criteria, and the third one is the verification of the practical performance of the presented algorithms, when applied to large sparse matrices from current research problems.

\appendix

\section{Effect of the parameter errors in Algorithm 1}
In this section, in order to check the effects of the parameter errors in Algorithm 1, we present some numerical examples.

\begin{figure}[htbp]
  \centering
  \includegraphics[width=0.95\linewidth]{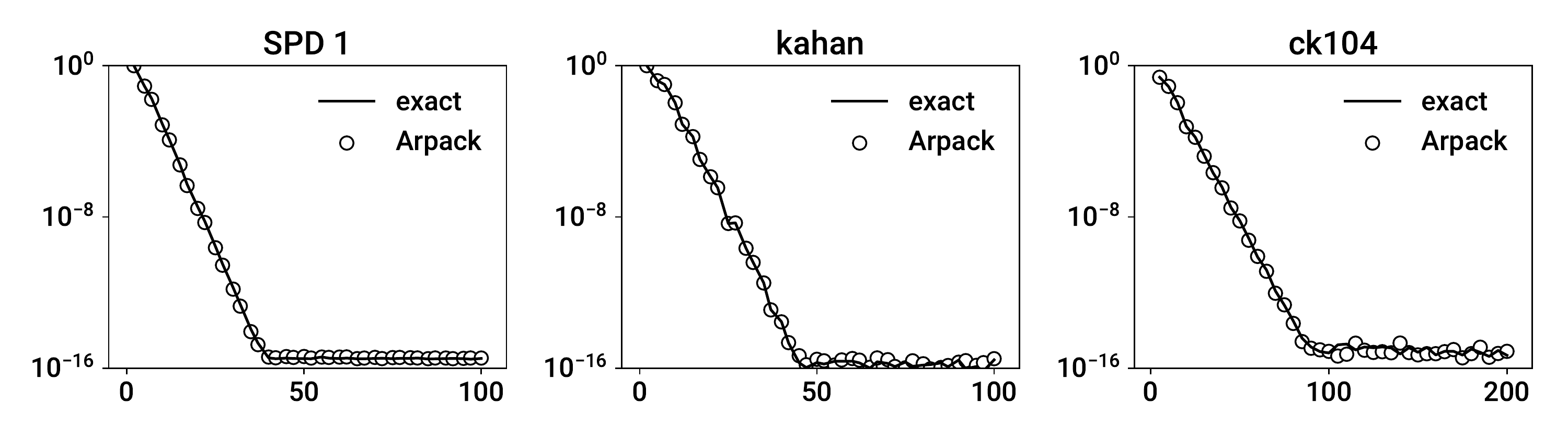}
  \caption{Convergence histories of the DE formula (Algorithm \ref{alg:1}) obtained using the exact estimations of $\|\tilde{A}-I\|_2,~\|\tilde{A}^{-1}\|_2,~\rho(\tilde{A})$ (using \texttt{eigvals} and \texttt{svdvals}) and obtained using \texttt{eigs} of \texttt{Arpack.jl}.
  The vertical axes show the relative error $\|\log(\tilde{A})-X\|_F/\|\log(\tilde{A})\|_F$, and the horizontal axes show the number of abscissas $m$.}
  \label{fig:2}
\end{figure}

The convergence histories of the DE formula are shown in Figure \ref{fig:2}.
Each graph shows two histories: one is obtained by using the \texttt{eigvals} and \texttt{svdvals} functions for computing $\rho(\tilde{A})$, $\|\tilde{A}-I\|_2$, and $\|\tilde{A}^{-1}\|_2$  as in Section 3;
the other one is obtained by using \texttt{Arpack.jl} with \texttt{tol=0.01}\footnote{
  The parameter \texttt{tol} defines the relative tolerance for convergence of Ritz values.
  In this examples, $\rho(\tilde{A})$ is obtained using \texttt{eigs(A, nev=1, which=:LM, tol=0.01)}, $\|\tilde{A}-I\|_2$ is obtained using \texttt{eigs((A-I)'*(A-I), nev=1, which=:LM, tol=0.01)}, and $\|\tilde{A}^{-1}\|_2$ is obtained using \texttt{eigs(A'*A, nev=1, which=:SM, tol=0.01)}, where \texttt{A} is the preconditioned matrix $\tilde{A}$.
}.
The figure shows that the behaviors of the histories are almost equal, and the effects of the errors did not appear.

In conclusion,
the parameters used at Step 2 of Algorithm \ref{alg:1} can be computed roughly.

\section{Adaptive quadrature algorithm based on the GL quadrature}
In this section we show an adaptive quadrature algorithm based on the GL quadrature for \eqref{eq:2-01}, using a technique from \cite{cardoso2012computation}.

Let $G(m):=\sum_{i=1}^m w_i F_{\mathrm{GL}}(u_i)$ be the results of the $m$-point GL quadrature, where $F_{\mathrm{GL}}(u_i)=[(1+u)(A-I)+2I]^{-1}$.
Using \cite[Eq. (4.6)]{cardoso2012computation}, the following error estimate can be applied:
\begin{align}\label{eq:5-2}
  \left\lVert G(2m)-\int_{-1}^1 F_{\mathrm{GL}}(t)\,\mathrm{d} t \right\rVert  \le  \left\lVert G(m) - G(2m) \right\rVert .
\end{align}

Based on \eqref{eq:5-2}, we present the algorithm designed to compute $\log(A)$ by automatically adding abscissas until the error is smaller than a given tolerance as Algorithm \ref{alg:3}.
\begin{algorithm}
  \caption{Computation $\log(A)$ based on Gauss--Legendre quadrature}
  \begin{algorithmic}[1]
    \Statex \textbf{Input:}  $A\in\mathbb{R}^{n\times n}$, The number of initial abscissas $m_0 \ge 2$, $\zeta>0$ for a tolerance of the truncation error.
    \Statex \textbf{Output:} $X\approx \log(A)$
    \State Set $F_{\mathrm{GL}}(u):=[(1+u)(A-I)+2I]^{-1}$
    \State Compute $\theta$, a lower bound on $\|\log(A)\|_2$.
    \State Compute abscissas $u_i$ and weights $w_i$ of $m_0$-point Gauss--Legendre quadrature ($i=1,\dots,m_0$)
    \State $G_0 = \sum_{i=1}^{m_0} w_i F_{\mathrm{GL}}(u_i)$
    \For{$k=0,1,2,\ldots$ until convergence}
      \State $m_{k+1}=2m_k$
      \State Compute abscissas $u_i$ and weights $w_i$ of $m_{k+1}$-point Gauss--Legendre quadrature ($i=1,\dots,m_{k+1}$)
      \State $G_{k+1} = \sum_{i=1}^{m_{k+1}} w_i F_{\mathrm{GL}}(u_i)$
      \If {$\|G_{k+1}-G_k\|/\theta \le \zeta$}
        \State $G=G_{k+1}$
        \State \textbf{break}
      \EndIf
    \EndFor
    \State $X = (A-I)G$
  \end{algorithmic}
  \label{alg:3}
\end{algorithm}

The computational cost of Algorithm \ref{alg:3} for dense $A$ is $[2(\sum_{i=0}^{k+1}m_i) + 2]n^3+\mathcal{O}(n^2) = [(2^{k+3}-2)m_0 + 2]n^3+\mathcal{O}(n^2)$, and the computational cost for $\log(A)\bm{b}$ with sparse $A$ is $(2^{k+2}-1)m_0c_{\mathrm{abscissa}} + c_{\mathrm{mul}} + c_{\mathrm{param}}$.
Under the assumption that the total computational cost largely depends on the coefficient of $c_{\mathrm{abscissa}}$, the computational cost of the Algorithm 2 will be smaller than that of Algorithm 3 when the convergence ratios of the DE formula and the GL quadrature are the same.

\end{document}